\documentclass[10pt]{article}
\usepackage{amssymb,amsmath,amsfonts,bbm,pifont,upgreek,bbold,accents} 
\usepackage[colorlinks=true]{hyperref}
\hypersetup{urlcolor=blue, citecolor=red}
 \usepackage[active]{srcltx}

\usepackage[dvips]{graphicx}
\font\teneufm=eufm10
\font\seveneufm=eufm7
\font\fiveeufm=eufm5
\newfam\eufmfam
\textfont\eufmfam=\teneufm
\scriptfont\eufmfam=\seveneufm
\scriptscriptfont\eufmfam=\fiveeufm
\def\eufm#1{{\fam\eufmfam\relax#1}}

\newcommand\beq[1]{ \begin{equation}\label{#1} }
\newcommand{\eeq}{ \end{equation} }

\newcommand\beqa[1]{ \begin{eqnarray} \label{#1}}
\newcommand{\eeqa}{ \end{eqnarray} }
\newcommand{\beqano}{ \begin{eqnarray*} }
\newcommand{\eeqano}{ \end{eqnarray*} }
\newcommand\arr[1]{\left\{\begin{array}{l}#1\end{array}\right.}
\renewcommand{\theequation}{\arabic{section}.\arabic{equation}}

\newtheorem{theorem}{Theorem}[section]
\newtheorem{definition}{Definition}[section]
\newtheorem{proposition}{Proposition}[section]
\newtheorem{lemma}{Lemma}[section]
\newtheorem{sublemma}{Sublemma}[section]
\newtheorem{remark}{Remark}[section]
\newtheorem{notationalremark}{Notational Remark}[section]
\newtheorem{corollary}{Corollary}[section]
\newtheorem{assumption}{Assumption}[section]
\newtheorem{claim}{Claim}[section]
\newtheorem{tools}{$\negsp\negsp$}[subsection]

\newcommand\thm[1]{ \begin{theorem}\label{#1}}
\newcommand\thmtwo[2]{ \begin{theorem}[#1]\label{#2}}
\newcommand\ethm{ \end{theorem} }
\newcommand\dfn[1]{ \begin{definition}\label{#1} \rm}
\newcommand\dfntwo[2]{ \begin{definition}[#1]\label{#2} \rm}
\newcommand\edfn{ \end{definition} }
\newcommand\pro[1]{ \begin{proposition}\label{#1}}
\newcommand\protwo[2]{ \begin{proposition}[#1]\label{#2}}
\newcommand\epro{ \end{proposition} }
\newcommand\lem[1]{ \begin{lemma}\label{#1}}
\newcommand\lemtwo[2]{ \begin{lemma}[#1]\label{#2}}
\newcommand\elem{ \end{lemma} }
\newcommand\sublem[1]{ \begin{sublemma}\label{#1}}
\newcommand\sublemtwo[2]{ \begin{sublemma}[#1]\label{#2}}
\newcommand\esublem{ \end{sublemma} }
\newcommand\rem[1]{ \begin{remark}\label{#1} \rm}
\newcommand\erem{ \end{remark} }
\newcommand\notrem[1]{ \begin{notationalremark}\label{#1} \rm}
\newcommand\enotrem{ \end{notationalremark} }
\newcommand\cor[1]{ \begin{corollary}\label{#1}}
\newcommand\cortwo[2]{ \begin{corollary}[#1]\label{#2}}
\newcommand\ecor{ \end{corollary} }
\newcommand\asmp[1]{ \begin{assumption}\label{#1}}
\newcommand\asmptwo[2]{ \begin{assumption}[#1]\label{#2}}
\newcommand\easmp{ \end{assumption} }
\newcommand\clm[1]{ \begin{claim}\label{#1}}
\newcommand\eclm{ \end{claim} }
\newcommand{\proof}{\par\medskip\noindent{\bf Proof\ }}
\newcommand\equ[1]{{\rm (\ref{#1})}}

\expandafter\chardef\csname pre amssym.def
at\endcsname=\the\catcode`\@
\catcode`\@=11
\def\undefine#1{\let#1\undefined}
\def\newsymbol#1#2#3#4#5{\let\next@\relax
 \ifnum#2=\@ne\let\next@\msafam@\else
 \ifnum#2=\tw@\let\next@\msbfam@\fi\fi
 \mathchardef#1="#3\next@#4#5}
\def\mathhexbox@#1#2#3{\relax
 \ifmmode\mathpalette{}{\m@th\mathchar"#1#2#3}%
 \else\leavevmode\hbox{$\m@th\mathchar"#1#2#3$}\fi}
\def\hexnumber@#1{\ifcase#1 0\or 1\or 2\or 3\or 4\or 5\or 6\or 7\or
8\or
 9\or A\or B\or C\or D\or E\or F\fi}
\ifcase\@ptsize
 \font\tenmsb=msbm10
 \font\sevenmsb=msbm7
 \font\fivemsb=msbm5
\or
 \font\tenmsb=msbm10 scaled \magstephalf
 \font\sevenmsb=msbm7 scaled \magstephalf
 \font\fivemsb=msbm5  scaled \magstephalf
\or
 \font\tenmsb=msbm10 scaled \magstep1
 \font\sevenmsb=msbm7 scaled \magstep1
 \font\fivemsb=msbm5 scaled \magstep1
\fi
\newfam\msbfam
\textfont\msbfam=\tenmsb
\scriptfont\msbfam=\sevenmsb
\scriptscriptfont\msbfam=\fivemsb
\edef\msbfam@{\hexnumber@\msbfam}
\def\Bbb#1{\fam\msbfam\relax#1}
\def\widehat#1{\setboxz@h{$\m@th#1$}%
 \ifdim\wdz@>\tw@ em\mathaccent"0\msbfam@5B{#1}%
 \else\mathaccent"0362{#1}\fi}
\def\widetilde#1{\setboxz@h{$\m@th#1$}%
 \ifdim\wdz@>\tw@ em\mathaccent"0\msbfam@5D{#1}%
 \else\mathaccent"0365{#1}\fi}

\def\RIfM@{\relax\ifmmode}
\def\nonmatherr@#1{\errmessage{\string#1\space allowed only in math mode}}
\def\Bbb{\RIfM@\expandafter\Bbb@\else
 \expandafter\nonmatherr@\expandafter\Bbb\fi}
\def\Bbb@#1{{\Bbb@@{#1}}}
\def\Bbb@@#1{\fam\msbfam\relax#1}
\def\setboxz@h{\setbox\z@\hbox}
\def\wdz@{\wd\z@}
\catcode`\@=\csname pre amssym.def at\endcsname
\newcommand{\eg}{{\tt e.g.\,}}

\newcommand{\Giu}{{\bigskip\noindent}}
\newcommand{\nl}{{\smallskip\noindent}}
\newcommand{\noi}{{\noindent}}

\newcommand{\qed}{\hskip.5truecm
\vrule width 1.7truemm height 3.5truemm depth 0.truemm
\par\Giu}

\newcommand{\negsp}{\hspace{-.09truecm}}  

\newcommand{\dst}{\displaystyle}
\newcommand\ovl[1]{ \overline {#1} }

\newcommand{\torus}{ {\Bbb T}   }
\renewcommand{\natural}{ {\Bbb N}   }
\newcommand{\real}{ {\Bbb R}   }
\newcommand{\integer}{ {\Bbb Z}   }
\newcommand{\complex}{ {\Bbb C}   }

\renewcommand{\d}{ {\delta}   }

\renewcommand{\l}{ {\lambda}   }
\renewcommand{\L}{ {\Lambda}   }
\newcommand{\m}{ {\mu}   }

\newcommand{\p}{ {\pi}   }

\renewcommand{\r}{ {\rho}   }
\newcommand{\s}{ {\sigma}   }

\renewcommand{\t}{ {\tau}   }
\newcommand{\f}{ {\varphi}   }

\renewcommand{\o}{ {\omega}   }

\newcommand{\ff}{ {\sl f} }
\newcommand{\cR}{ {\cal R} }

\newcommand{\cL}{ {\cal L} }

\newcommand{\cN}{ {\cal N} }
\newcommand{\cP}{ {\cal P} }
\newcommand{\cI}{ {\cal I} }
\newcommand{{\cJ}}{ {\cal J} }

\newcommand{\cX}{ \Xi }

\newcommand\ppu{{ (1) }}

\newcommand\ppj{{ (j) }}

\newcommand\ul{{\uplambda}}

\newcommand\up{{\rm p}}
\newcommand\uq{{\rm q}}

\newcommand\cO{{\cal O}}

\newcommand\HH{{\rm H}}
\newcommand\II{{\rm I}}

\newcommand\hh{{\rm h}}

\newcommand\JJ{{\rm J}}
\newcommand\rr{{\rm r}}

\newcommand\ii{{\rm i}}

\newcommand\EE{{\rm E}}

\begin{document}

\title{A normal form without small divisors\\
(Draft)}

\author{  
Gabriella Pinzari\\\\
\footnotesize{Dipartimento di Matematica Universit\`a di Padova}\\
\footnotesize{gabriella.pinzari@math.unipd.it}
}

\date{April 2, 2017}
\maketitle
\vskip.1in
\noi


\maketitle
\tableofcontents
\abstract{Following the techniques of \cite{poschel93}, we formulate a Normal Form Lemma suited to close-to-be-integrable Hamiltonian systems where not all the coordinates are action--angles. The Lemma turns to be useful in the theory of KAM tori of Sun-Earth-Asteroids systems (work in progress of the author; arXiv: 1702.03680).}

\section{
Set Up
}
\setcounter{equation}{0}
\renewcommand{\theequation}{\arabic{equation}}
Consider the $2(n+m+1)$--dimensional phase space $$\cP=\cR\times \cI \times \cX\times \torus^n\times B^{2m}_{\d}$$ 
where $\cR\subset\real$, $\cI\subset \real^n$, $0\in\Xi\subset\real$ are open, connected and bounded, while $B^{2m}_{\d}$ denotes the ball of radius ${\d}$ in $\real^{2m}$ centered at $0\in \real^{2m}$. Le $\cP$ be 
equipped with set of canonical coordinates $(\rr,\II,x,\f,p,q)\in \cP$ with respect to the standard two--form
$$\Omega=d\rr\wedge dx+d\II\wedge d\f+dp\wedge dq=d\rr\wedge dx+\sum_{i=1}^nd\II_i\wedge d\f_i+\sum_{j=1}^mdp_j\wedge dq_j$$
and
 consider, on $\cP$, a Hamiltonian of the form

\beq{full system}\HH(\rr,\II,p, x,\f,q)=\HH_0(\rr,\II,\JJ(p,q))+f(\rr,\II,p, x,\f,q)\eeq
where
$$\JJ(p,q)=(p_1q_1,\cdots, p_mq_m)\ .$$
Note that {\it we are not assuming that $f$ is periodic in $x$}.

\nl
Setting $f$ to zero, the Hamiltonian $\HH=\HH_0$ has the  motions 
\beq{unperturbed motion}\arr{\dst
\II=\II_\star\\
\dst  \rr=\rr_\star\\
\dst p=p_\star e^{-\omega_\JJ(t-t_\star)}}\qquad \arr{
\dst x=x_\star+\o_r (t-t_\star)\\
\dst  \f=\f_\star+\o_\II (t-t_\star)\\
\dst q=q_\star e^{\omega_\JJ(t-t_\star)}\\
}\eeq
where
$$\o_{\rr,\II,\JJ}:=\partial_{\rr,\II,\JJ}\hh(\rr,\II,\JJ)$$

\nl
We consider the problem of the continuation of such motions to the full system \equ{full system}.

\nl
The problem may be regarded as a generalization of problems that have been widely investigated in the framework of {\sc kam} and Nekhorossev theory.

\nl
In fact, if $\HH_0$ was taken to be independent of $\rr$, 
 we would be  in the setting of (partially hyperbolic) {\sc kam}
theory, where the perturbing function will depend, in addition, on the ``degenerate'' couple $(\rr, x)$. Such case has been investigated in the literature, starting with V.I. Arnold and N.N. Nekhorossev  \cite{arnold63, nehorosev77}.
 Refinements have been given by L. Chierchia and G. Pinzari in the case of properly--degenerate {\sc kam} theory \cite{chierchiaPi10}, by J. P\"oschel in the case of  Nekhorossev theory \cite{poschel93}.  Such papers are addressed to the study of
 Hamiltonian systems (named ``properly--degenerate'') of the form
$$\HH=\HH_0(\II)+f(\II, \f,\up,\uq)\qquad (\II, \f,\up,\uq)\in \cI\times \torus^n\times B^{2m}_{\d} $$
where the unperturbed part $\HH_0$ has strictly less degrees of freedom than the whole system. For such systems standard techniques do not apply since, on one side, as for  {\sc kam} theory,  usual non-degeneracy assumptions are strongly prevented and, on the other site, as for Nekhorossev theory, one has to control the variation of the ``degenerate'' coordinates'' $(\up,\uq)$.  For the way how such difficulties have been overcome, we refer to the dedicated literature (recalled in \cite{chierchiaPi10} and references therein).

\nl
The generalization studied in this paper with respect to the previous mentioned  cases is precisely related to the r\^ole of the coordinate $x$:  we are {\it not} assuming that this is a periodic coordinate, henceforth, standard  {\sc kam} theories do not apply. In this setting, one cannot reasonably expect, at least in general, that its linear motion of $x$ in \equ{unperturbed motion} is preserved {\it at any time}. 

\nl
As an example, let us look at the clock Hamiltonian 
$$\HH_\varepsilon=\frac{r^2}{2}+\varepsilon^2 \frac{x^2}{2}\ .$$
For $\varepsilon=0$, $\HH_\varepsilon$ reduces to the free hamiltonian
$\dst\HH_0=\frac{\rr^2}{2}$
whose motions are
$$r_0(t)=r_\star\qquad x_0(t)=x_\star+r_\star(t-t_\star)\ .$$
However, when $\varepsilon\ne 0$, the motions of  $\HH_\varepsilon$, given by
$$r_\varepsilon(t)=r_\star\cos\varepsilon(t-t_\star)-\varepsilon x_\star\sin\varepsilon(t-t_\star)\qquad x_\varepsilon(t)=\frac{r_\star}{\varepsilon}\sin\varepsilon(t-t_\star)+ x_\star\cos\varepsilon(t-t_\star)$$
are effectively close one to the one of $\HH_0$  for $|t-t_\star|$ of the order $\varepsilon^{-1}$.
For larger times, the two Hamiltonians generate a completely different dynamics, since the former has only unbounded motions, while the latter has bounded ones. The same conclusion could be reached, instead of solving the motion equations,  looking at the  phase portrait of $\HH_\varepsilon$, which consists of ellipses with semi--axes $\sqrt{2\EE}$, $\varepsilon^{-1}\sqrt{2\EE}$ which tend to the straight lines $\rr=\sqrt{2\EE}$ as $\varepsilon\to 0$.

\nl
Similarly, one sees that still in the case of the Hamiltonian
$$\widetilde\HH_\varepsilon=\frac{r^2}{2}-\varepsilon^2 \frac{x^2}{2}$$
which has unbounded motions for all $\varepsilon$, the dynamics of $\widetilde\HH_0$ and $\widetilde\HH_\varepsilon$ with $\varepsilon\ne 0$ are very far one from the other for $|t-t_\star|\gg \varepsilon^{-1}$.

\nl
For these reasons, we divide the problem of the study of the dynamics of the full Hamiltonian \equ{full system} in two steps. As a first step, which is actually the purpose of this note, we consider the {\it intermediate problem} of constructing a {\it normal form} for $\HH$ for  very large (exponentially long) times, without any attempt to normalize the evolution of the couple $(r, x)$. Such normal form will be defined on a suitable sub--domain $\cP_N\subset\cP$, and will be of the kind
$$\HH_{N}=\HH_0(\rr,\II, pq)+\HH_1(\rr,\II, pq, x)+f_{N}(\rr, \II, p, x,\f,q)$$
where $f_{N}$ is a very (exponentially) small remainder which we shall quantify. Clearly, when dealing with concrete applications,  such step should be followed by a second step where one verifies that the evolution generated by $\HH_N$ remains in the prescribed domain $\cP_{N}$ for all such time. Such idea of ``a posteriori check'' goes back to N.N. Nekhorossev \cite{nehorosev77}, who indeed was able to establish
the validity, to the ${\rm N}$--body problem Hamiltonian written in Poincar\'e coordinates
$$\HH_\cP=\hh_0(\L)+\ff_\cP(\L,\ul,\up,\uq)$$
 of a normal form of the kind
$$\HH_N=\hh_0(\L)+\hh_1(\L,\up,\uq)+{\rm f}_N(\L,\ul,\up,\uq)$$
where ${\rm f}_N$ is exponentially small,
just controlling that the  ``degenerate'' coordinates $(\up,\uq)$  did not escape their domain for all that time. Before stating our result, let us fix the following

\paragraph{Notations}

We consider the complex neighborhood 
$$\cP_{r,\r,\xi,s,\d}=\cR_r\times \cI_\r \times \Xi_\xi\times \torus^n_s\times B^{2m}_{\d}\ ,$$
of $\cP$ where, as usual, $A_\theta:=\cup_{x_0\in A}\{B_\theta(x_0)\}$, while $\torus_s:=\torus+\ii [-s,s]$, with $\torus:=\real/(2\p\integer)$ the standard torus. 

\nl
We denote as ${\cal O}_{r,\r,\xi,s,\d}$ the set of complex holomorphic functions $\phi:\ \cP_{\hat r,\hat \r,\hat \xi,\hat s,\hat d}\to \complex$ for some $\hat r>r$, $\hat\r>\r$, $\hat\xi>\xi$, $\hat s>s$, $\hat\d>\d$.

\nl
We equip  ${\cal O}_{r,\r,\xi,s,\d}$  with the norm
$$\|\phi\|_{r,\r,\xi,s,\d}:=\sum_{k,h,j}\|\phi_{khj}\|_{r,\r,\xi}e^{s|k|}\d^{h+j}$$
where $\phi_{khj}(\rr,\II, x)$ are the coefficients of the Taylor--Fourier expansion
$$\phi=\sum_{k,h,j}\phi_{khj}(\rr,\II,x)e^{\ii k s}p^h q^j\ ,$$
and $\|\phi_{khj}\|_{r,\r,\xi}:=\sup_{\cR_r\times \cI_\r\times \Xi_\xi}|\phi_{khj}|$. Observe that  $\|g_{khj}\|_{r,\r,\xi}$ is well defined because of the boundedness of $\cR$, $\cI$ and $\Xi$, while $\|\phi\|_{r,\r,\xi,s,\d}$ is well defined by the usual properties of  holomorphic functions.

\nl
For a given vector--valued  function $\underline\phi=(\phi_1,\cdots, \phi_k)\in {\cal O}_{r,\r,\xi,s,\d}^k$, we let
$$\|\underline\phi\|_{r,\r,\xi,s,\d}:=\sum_{i=1}^k \|\phi_i\|_{r,\r,\xi,s,\d}\ .$$

\nl
If $\phi\in {\cal O}_{r,\r,\xi,s,\d}$, we define its ``off--average''  and   ``average''  as
$$\widetilde\phi:=\sum_{k,h,j:\atop (k,h-j)\ne (0,0)}g_{khj}(\rr,\II, x)e^{\ii k s}p^h q^j\ ,\qquad
\ovl\phi:=\phi-\widetilde\phi
\ .$$
Then we define the ``zero--average'' and the the ``normal'' classes
as
\beqa{zero average}&&{\cal Z}_{r,\r,\xi,s,\d}:=\{\phi\in {\cal O}_{r,\r,\xi,s,\d}:\quad \phi=\widetilde\phi\}=\{\phi\in {\cal O}_{r,\r,\xi,s,\d}:\quad \ovl\phi=0\}\\
\label{phi independent}&&{\cal N}_{r,\r,\xi,s,\d}:=\{\phi\in {\cal O}_{r,\r,\xi,s,\d}:\quad\phi=\ovl\phi\}=\{\phi\in {\cal O}_{r,\r,\xi,s,\d}:\quad \widetilde\phi=0\}\ .\eeqa
respectively. Obviously, one has the decomposition
$${\cal O}_{r,\r,\xi,s,\d}={\cal Z}_{r,\r,\xi,s,\d}\oplus {\cal N}_{r,\r,\xi,s,\d}\ .$$

\paragraph{Result} We  assume we are given a Hamiltonian system of the form \equ{full system}, 
where $\HH_0\in {\cal N}_{r,\r,\xi,s,\d}$, is $x$--independent and $f\in {\cal O}_{r,\r,\xi,s,\d}$.

\nl
We shall prove the following result.

\begin{lemma}[Normal Form Lemma]\label{normal form lemma}
There exists a number ${\tt c}_{n,m}\ge 1$ such that, for any $N\in \natural$ such that the following inequalities are satisfied
\beq{normal form assumptions} 4N{\cal X}\|\frac{\o_\II}{\o_\rr}\|_{r,\r}<
s
\ ,\quad
4N{\cal X}\|\frac{\o_\JJ}{\o_\rr}\|_{r,\r}<
1\ ,\quad {\tt c}_{n,m}N\frac{{\cal X}}{{\eufm d}}   \|\frac{1}{\o_\rr}\|_{r,\r} \| f\|_{r,\r,\xi,s,\d} <1 \eeq
with ${\eufm d}:=\min\big\{\r\s, r\xi, {\d}^2\big\}$, ${\cal X}:=\sup\big\{|x|:\ x\in \Xi_\xi\big\}$ and $\omega_{\rr,\II,\JJ}:=\partial_{\rr,\II,\JJ} \HH_0$, one can find an operator $$\Psi_N:\quad {\cal O}_{r,\r,\xi,s,\d}\to \cO_{1/3 (r, \r,\xi, s, \d)}$$
which carries $\HH$ to
$$\HH_N:=\Psi_N[\HH]=\HH_0+\HH_1+f_N$$
where
$\HH_1\in \cN_{1/3 (r, \r,\xi, s, \d)}$, $f_N\in \cO_{1/3 (r, \r,\xi, s, \d)}$ and, moreover, the following inequalities hold 
\beqa{thesis}
&&\|\HH_1-\ovl f\|_{1/3 (r, \r,\xi, s, \d)}\le {\tt c}_{n,m} \frac{{\cal X}}{\eufm d}\|\frac{1}{\o_\rr}\|_{r,\r} \|\widetilde f\|_{r,\r,\xi,s,\d}\| f\|_{r,\r,\xi,s,\d}\nonumber\\
&&    \|f_N\|_{1/3 (r, \r,\xi, s, \d)}\le \frac{1}{2^{N+1}} \|f\|_{r,\r,\xi,s,\d}\ .\eeqa

\end{lemma}

\nl
 The main point of Lemma \ref{normal form lemma} is that it holds {\it without small denominators}. However, we set an additional requirement that  the frequencies $\omega_\II$ and $\omega_\JJ$  are {\it small} (compare the two former inequalities in \equ{normal form assumptions}). Such assumption, that may seem too restrictive in general, has however many chances of being satisfied in the case of system arising from Celestial Mechanics, because, due to the proper degeneracy recalled above, very often, one has to deal with an ``effective system'' whose  unperturbed part includes  some manipulation of the perturbing function, which is naturally small.  We shall show a situation where indeed this is the case in a forthcoming paper.
 
 \nl
We now aim to give an account of the basic idea that enabled us to avoid the  small--divisor problem, by underlying the differences with the ``standard'' situation. We call so the situation, largely studied in the aforementioned papers, where $f$ is $x$--periodic, and one looks for a $\Psi_N$ which is also $x$--periodic.

\nl The beginning is just as in the standard case. We follow  the well--settled framework acknowledged to J\"urgen P\"oschel \cite{poschel93}. As in \cite{poschel93}, we shall obtain Lemma \ref{normal form lemma}  via iterate applications of one--step transformations (Iterative Lemma, see below) where the dependence of $\f$ and $(p,q)$ other than the combinations  $\JJ(p,q)$ is eliminated at higher and higher orders. It  goes as follows.

\nl
We assume that, at a certain step, we have a system of the form
\beq{step i}\HH=\HH_0(\rr,\II, \JJ(p,q))+g(\rr,\II, \JJ(p,q), x)+f(\rr,\II,  x, \f, p,q)\eeq
where  $f\in {\cal O}_{r,\r,\xi,s,\d}$, while $\HH_0$, $g\in {\cal N}_{r,\r,\xi,s,\d}$, with $\HH_0$ is independent of $x$ (the first step corresponds to take $g\equiv0$). 

\nl
After splitting $f$ on its Taylor--Fourier basis 
$$f=\sum_{k,h,j} f_{khj}(\rr,\II,x)e^{\ii k \f}p^h  q^j\ .$$
one looks for a time--1 map 
$$\Phi=e^{\cL_\phi}$$
generated by a small Hamiltonian 
$\phi$
which will be taken  in the  class ${\cal Z}_{r,\r,\xi,s,\d}$ in \equ{zero average}.
One lets
\beq{exp}\phi=\sum_{(k,h,j):\atop{(k,h-j)\ne (0,0)}} \phi_{khj}(\rr,\II,x)e^{\ii k \f}p^h q^j\ .\eeq

\nl
The operation
$$\phi\to \{\phi,\HH_0\}$$
acts diagonally  on the monomials  in the expansion \equ{exp}, carrying
\beq{diagonal}\phi_{khj}\to -\big(\o_\rr\partial_x \phi_{khj}+\l_{khj} \phi_{khj}\big)\ ,\quad {\rm with}\quad \l_{khj}:=(h-j)\cdot\o_\JJ+\ii k\cdot\o_\II\ .\eeq
Therefore, one defines
$$\{\phi,\HH_0\}=:-D_\omega\phi\ .$$
The formal application of $\Phi=e^{\cL_\phi}$ yields:
\beqa{f1}
e^{\cL_\phi} \HH&=&e^{\cL_\phi} (\HH_0+g+f)=\HH_0+g-D_\omega \phi+f+\Phi_2(\HH_0)+\Phi_1(g)+\Phi_1(f)
\eeqa
where the $\Phi_h$'s are the queues of $e^{\cL_\phi}$, defined in Section \ref{Time--one flows}.

\nl
Next, one requires that   the residual term $-D_\omega \phi+f$ lies in the class ${\cal N}_{r,\r,\xi,s,\d}$ in \equ{phi independent}. This
amounts  to solve 
the ``homological'' equation
\beq{homological equation}\widetilde{(-D_\omega \phi+f\big)}=0\eeq
for $\phi$.

\nl
Since we have chosen $\phi\in{\cal Z}_{r,\r,\xi,s,\d}$, by \equ{diagonal}, we have that also $D_\omega \phi\in{\cal Z}_{r,\r,\xi,s,\d}$. So, Equation \equ{homological equation} becomes
\beq{homol}-D_\omega \phi+\widetilde f=0\ .\eeq
In terms of  the Taylor--Fourier modes, the equation becomes
\beq{homol eq} \o_\rr\partial_x \phi_{khj}+\l_{khj} \phi_{khj}=f_{khj}\qquad \forall\ (k,h,j):\ (k,h-j)\ne (0,0)\ .\eeq

\nl In the standard situation, one typically proceeds to solve such equation via Fourier series:
\beq{periodic case}f_{khj}(\rr,\II,x)=\sum_{\ell}f_{khj\ell}(\rr,\II)e^{\ii \ell x}\ ,\qquad \phi_{khj}(\rr,\II,x)=\sum_{\ell}\phi_{khj\ell}(\rr,\II)e^{\ii \ell x}\eeq
so as to find
$\dst\phi_{khj\ell}=\frac{f_{khj\ell}}{\m_{khj\ell}}$
with the usual  denominators $\m_{khj\ell}:=\l_{khj}+\ii\ell \omega_\rr$ which one requires not to vanish 
via, \eg, a ``diophantine inequality'' to be held for all $(k,h,j,\ell)$ with $(k,h-j)\ne (0,0)$. Observe that, in the classical case, there is not much freedom in the choice of $\phi$. In fact, such solution is determined up to solutions  of the homogenous equation
\beq{homogeneous}D_\omega\phi_0=0\eeq
which, in view of the Diophantine condition, has the only trivial solution $\phi_0\equiv0$.

\nl
{\it The situation is different if $f$ is not periodic in $x$, or $\phi$ is not needed so}. In such a case, it is possible to find a solution of \equ{homol eq}, corresponding to a non--trivial solution of \equ{homogeneous},  where small divisors do not appear.

\nl
This is
\beq{solution}\phi_{khj}(\rr,\II, x)=
\frac{1}{\o_\rr}\int_0^xf_{khj}(\rr,\II,\t)e^{\frac{\l_{khj}}{\o_\rr}(\t-x)}d\t\qquad \forall\ (k,h,j):\ (k,h-j)\ne (0,0)
\eeq
and $\phi_{0hh}(\rr,\II, x)\equiv0$. Note that in the particular case that $f$ is periodic in $x$, and hence it affords an expansion like \equ{periodic case}, the solution \equ{solution} may be  written as
$$\phi_{khj}=e^{-\frac{\l_{khj}}{\o_\rr}x}\sum_{\ell}f_{khj\ell}(\rr,\II) \frac{e^{\frac{\m_{khj\ell}}{\o_\rr}x}-1}{\m_{khj\ell}}=:e^{-\frac{\l_{khj}}{\o_\rr}x}\widehat\phi_{khj}(\rr,\II,x, p,q)$$
As expected, such a solution provides, via \equ{exp}, a function $\phi$ that, in general,  is {\it not} periodic in $x$ for all $(\rr,\II,\f, p,q)$ in their domain. Indeed, 
{\it under the genericity assumption that 
the $\widehat\phi_{khj}$'s have no other common zero than $x=0$,}
since such  $\widehat\phi_{khj}$'s are periodic in $x$, we have that the  $\phi_{khj}$'s are so only for $(\rr,\II,p,q)$ such that $\frac{\l_{khj}}{\o_\rr}\in \ii\integer$. Henceforth, $\phi$ is $x$--periodic
only on the subset
$(\rr,\II,p,q, \f)\in \cR_{es}\times \torus^n$, where $\cR_{es}$ is the zero--measure subset of $\cR\times\cI\times B^{2m}_\d$ where $\frac{\l_{khj}}{\o_\rr}\in \ii\integer$ for all $(k,j,h)$ such that $(k,h-j)\ne (0,0)$. 

\nl
We conclude  with a comment on the necessity of the two first inequalities in \equ{normal form assumptions}:  the formula \equ{solution} involves some loss of analyticity for $\phi$ whose strength we will evaluate to be of the order of the maximum of ${{\cal X}\|\frac{\o_\II}{\o_\rr}\|_{r,\r}}$, ${{\cal X}\|\frac{\o_\JJ}{\o_\rr}\|_{r,\r}}$.

\newpage
\section{Proofs}\label{Time--one flows}
\begin{definition}[Time--one flows and their queues]\rm
Let $\cL_\phi(\cdot):=\big\{\phi, \cdot\big\}$, where $\{f, g\}:=\sum_{i=1^k}(\partial_{p_i}f\partial_{q_i}g-\partial_{p_i}g\partial_{q_i}f)$, where $\Omega=\sum_{i=1}^k dp_i\wedge dq_i$ is the standard two--form, denotes Poisson parentheses.

\nl
For a given $\phi\in {\cal O}_{r,\r,\xi,s,\d}$, we denote as $\Phi_h$, $\Phi$ the formal series

\beq{queue}\Phi_h:=\sum_{j\ge h}\frac{{\cal L}_\phi^j}{j!}\ \qquad \Phi:=\Phi_0\ .\eeq
It is customary to let, also $\Phi:=e^{\cL_\phi}$.
\end{definition}

\begin{lemma}[\cite{poschel93}]\label{base lemma}
There exist an integer number $\ovl{\tt c}_{n,m}$ such that, for any $\phi\in {\cal O}_{r,\r,\xi,s,\d}$ and any $r'<r$, $s'<s$, $\r'<\r$, $\xi'<\xi$, $\d'<\d$ such that
$$\frac{\ovl{\tt c}_{n,m}\|\phi\|_{r,\r,\xi,s,\d}}{d}<1\qquad d:=\min\big\{\r'\s', r'\xi', {\d'}^2\big\}$$
then  the series in \equ{queue} converge uniformly so as to define 
the family $\{\Phi_h\}_{h=0,1,\cdots}$ of operators 
 $$\Phi_h:\quad  {\cal O}_{r,\r,\xi,s,\d}\to \cO_{r-r',\r-\r',\xi-\xi',s-s',\d-\d'}\ .$$ 
 Moreover, the following bound holds (showing, in particular, uniform  convergence):

\beq{geometric series}\|\cL^j_\phi[g]\|_{r-r',\r-\r',\xi-\xi',s-s',\d-\d'}\le j!\big(\frac{\ovl{\tt c}_{n,m}\|\phi\|_{r,\r,\xi,s,\d}}{d}\big)^j\|g\|_{r,\r,\xi,s,\d}\ .\eeq
for all $g\in{\cal O}_{r,\r,\xi,s,\d}$.
\end{lemma}

\begin{remark}[\cite{poschel93}]\rm The bound \equ{geometric series} immediately implies
\beq{h power}\|\Phi_hg\|_{r-r',\r-\r',\xi-\xi',s-s',\d-\d'}\le \frac{\big(\frac{c\|\phi\|_{r,\r,\xi,s,\d}}{d}\big)^h}{1-\frac{c\|\phi\|_{r,\r,\xi,s,\d}}{d}}\|g\|_{r,\r,\xi,s,\d}\qquad \forall g\in {\cal O}_{r,\r,\xi,s,\d}\ .\eeq
\end{remark}

\begin{lemma}[Iterative Lemma]\label{iterative lemma}
There exists a number $\widetilde{\tt c}_{n,m}>1$ such that the following holds. For any choice  of positive numbers  $r'$, $\r'$, $s'$, $\xi'$. $\d'$ satisfying
\beqa{ineq1}
&&{2r'<r\ ,\quad 2\r'<\r\ ,\quad 2\xi'<\xi}\\
\label{ineq2}
&& {2s'<s\ ,\quad 2\d'<\d\ ,\quad {\cal X}\|\frac{\o_\II}{\o_\rr}\|_{r,\r}<s-2s'\ ,\quad
{\cal X}\|\frac{\o_\JJ}{\o_\rr}\|_{r,\r}<
\log\frac{\d}{2\d'} }
\eeqa
and and provided that the  following inequality holds true
\beqa{smallness}
 \widetilde{\tt c}_{n,m}\frac{{\cal X}}{d}    \|\frac{1}{\o_\rr}\|_{r,\r} \|\widetilde f\|_{r,\r,\xi,s,\d} <1\qquad d:=\min\big\{\r'\s', r'\xi', {\d'}^2\big\}
\eeqa
 one can find an operator $$\Phi:\quad {\cal O}_{r,\r,\xi,s,\d}\to \cO_{r_+,\r_+,\xi_+,s_+,\d_+}$$
 with
 $${r_+:=r-2r'\ ,\quad \r_+:=\r-2\r'\ ,\quad \xi_+:=\xi-2\xi'\ ,\quad s_+:=s-2s'-{\cal X}\|\frac{\o_\II}{\o_\rr}\|_{r,\r}\ ,\quad \d_+:=\d e^{-{\cal X}\|\frac{\o_\JJ}{\o_\rr}\|_{r,\r}}-2\d'}$$
which carries the Hamiltonian  $\HH$ in \equ{step i} to
$$\HH_+:=\Phi[\HH]=\HH_0+g+\ovl f+f_+$$
where
\beq{bound}\|f_+\|_{r_+,\r_+, \xi_+, s_+,\d_+}\le \widetilde{\tt c}_{n,m}\frac{{\cal X}}{d}    \|\frac{1}{\o_\rr}\|_{r,\r} \|\widetilde f\|_{r,\r,\xi,s,\d}\| f\|_{r,\r,\xi,s,\d}+ \|\{\phi, g\}\|_{r_1-r',\r_1-\r',\xi_1-\xi',s_1-s',\d_1-\d' }
\eeq
with $$r_1:=r\ ,\quad \r_1:=\r\ ,\quad \xi_1:=\xi\ ,\quad s_1:=s-{\cal X}\|\frac{\o_\II}{\o_\rr}\|_{r,\r}\ ,\quad \d_1:=\d e^{-{\cal X}\|\frac{\o_\JJ}{\o_\rr}\|_{r,\r}}$$ for a suitable $\phi\in \cO_{r_1,\r_1,\xi_1,s_1,\d_1}$ verifying \beq{bound on phi}\|\phi\|_{r_1,\r_1,\xi_1,s_1,\d_1 }\le \frac{{\cal X}}{d}    \|\frac{1}{\o_\rr}\|_{r,\r} \|\widetilde f\|_{r,\r,\xi,s,\d}\ .\eeq
\end{lemma}

\proof 
Let $\ovl{\tt c}_{n,m}$ be as in Lemma \ref{base lemma}. We shall choose $\widetilde{\tt c}_{n,m}$ suitably large with respect to $\ovl{\tt c}_{n,m}$.

\nl
Let $\phi_{khj}$ as in \equ{solution}. 
Let us fix
\beq{ovl param}0<\ovl r\le r\ ,\quad 0<\ovl\r\le \r\ ,\quad 0<\ovl\xi\le \xi\ ,\quad 0<\ovl s< s\ ,\quad 0<\ovl\d< \d\eeq
and assume that
\beq{ovl param1}{\cal X}\|\frac{\o_\II}{\o_\rr}\|_{ r,\r} \le s-\ovl s\ ,\qquad {\cal X}\|\frac{\o_\JJ}{\o_\rr}\|_{ r,\r} \le \log\frac{\d}{\ovl \d}\ .\eeq
Then we have
$$\|\phi_{khj}\|_{\ovl r,\ovl\r,\ovl\xi}\le  \|\frac{1}{\o_\rr}\|_{\ovl r,\ovl\r} \|f_{khj}\|_{\ovl r,\ovl\r,\ovl\xi}\|\int_0^x |e^{-\frac{\l_{khj}}{\o_\rr}\t}|\|_{\ovl r,\ovl\r,\ovl\xi}d\t\le{\cal X}   \|\frac{1}{\o_\rr}\|_{\ovl r,\ovl\r}  \|f_{khj}\|_{\ovl r,\ovl\r,\ovl\xi} e^{{\cal X}\|\frac{\l_{khj}}{\o_\rr}\|_{\ovl r,\ovl\r}}\ .$$
Since
$$\|\frac{\l_{khj}}{\o_\rr}\|_{\ovl r,\ovl\r} \le (h+j)\|\frac{\o_\JJ}{\o_\rr}\|_{\ovl r,\ovl\r} +|k|\|\frac{\o_\II}{\o_\rr}\|_{\ovl r,\ovl\r} 
$$
we have, definitely,
$$\|\phi_{khj}\|_{\ovl r,\ovl\r,\ovl\xi}\le {\cal X}   \|\frac{1}{\o_\rr}\|_{\ovl r,\ovl\r}  \|\widetilde f_{khj}\|_{\ovl r,\ovl\r,\ovl\xi}
e^{(h+j){\cal X}\|\frac{\o_\JJ}{\o_\rr}\|_{\ovl r,\ovl\r} +|k|{\cal X}\|\frac{\o_\II}{\o_\rr}\|_{\ovl r,\ovl\r} }\ .
$$
which yields (after multiplying by $e^{|k|\ovl s}(\ovl\d)^{j+h}$ and summing over $k$, $j$, $h$ with $(k,h-k)\ne (0,0)$) to
$$\|\phi\|_{\ovl r,\ovl\r,\ovl\xi, \ovl s,\ovl\d}\le{\cal X}  \|\frac{1}{\o_\rr}\|_{\ovl r,\ovl \r,\ovl \xi} \|\widetilde f\|_{\ovl r,\ovl \r,\ovl \xi,\ovl s+{\cal X}\|\frac{\o_\II}{\o_\rr}\|_{\ovl r,\ovl\r} ,\ovl \d e^{{\cal X}\|\frac{\o_\JJ}{\o_\rr}\|_{\ovl r,\ovl\r} }}\ . $$
Note that the right hand side is well defined because of \equ{ovl param1}.
In the case of the choice
$$\ovl r=r=:r_1\ ,\quad\ovl\r=\r=: \r_1\ ,\quad \ovl\xi=\xi=:\xi_1\ ,\quad \ovl s=s-{\cal X}\|\frac{\o_\II}{\o_\rr}\|_{r,\r}=:s_1\qquad \ovl\d=\d e^{-{\cal X}\|\frac{\o_\JJ}{\o_\rr}\|_{r,\r}}=:\d_1$$
(which, in view of the two latter inequalities in  \equ{ineq2}, satisfies \equ{ovl param}--\equ{ovl param1}) the inequality becomes \equ{bound on phi}.
An application of Lemma \ref{base lemma},with $r$, $\r$, $\xi$, $s$, $\d$ replaced by $r_1-r'$, $\r_1-\r'$, $\xi_1-\xi'$, $s_1-s'$, $\d_1-\d'$,  concludes  with a suitable choice of $\widetilde{\tt c}_{n,m}>\ovl{\tt c}_{n,m}$ and (by \equ{f1}) $$f_+:=\Phi_2(\HH_0)+\Phi_1(g)+\Phi_1(f)\ .$$
Observe that the bound \equ{bound} follows from  Equations \equ{h power},  \equ{geometric series} and the identities
$$\Phi_2[\HH_0]=\sum_{j=2}^\infty \frac{\cL^j_\phi(\HH_0)}{j!}=\sum_{j=1}^\infty \frac{\cL^{j+1}_\phi(\HH_0)}{(j+1)!}=-\sum_{j=1}^\infty \frac{\cL^{j}_\phi(\widetilde f)}{(j+1)!}$$
$$\Phi_1[g]=\sum_{j=1}^\infty \frac{\cL^j_\phi(g)}{j!}=\sum_{j=0}^\infty \frac{\cL^{j+1}_\phi(g)}{(j+1)!}=-\sum_{j=0}^\infty \frac{\cL^{j}_\phi(g_1)}{(j+1)!}$$
with $g_1:=\cL_\phi(g)=\{\phi, g\}$.\qed

\nl
The proof of Lemma \ref{normal form lemma} goes through iterate applications of Lemma \ref{iterative lemma}. At this respect, we premise the following

\begin{remark}\label{stronger iterative lemma}\rm
Replacing conditions in \equ{ineq2} with the stronger ones

\beq{new cond}
{3s'<s\ ,\quad 3\d'<\d\ ,\quad {\cal X}\|\frac{\o_\II}{\o_\rr}\|_{r,\r}<s'\ ,\quad
{\cal X}\|\frac{\o_\JJ}{\o_\rr}\|_{r,\r}<
\frac{\d'}{\d}} \eeq

\nl
(and keeping \equ{ineq1}, \equ{smallness} unvaried)  one can take, for $s_+$, $\d_+$,  $s_1$, $\d_1$ the simpler expressions

$${s_{+\rm new}=s-3s'\ ,\quad \d_{+\rm new}=\d-3\d'\ ,\quad s_{1\rm new}:=s-s'\ ,\quad \d_{1\rm new}=\d-\d'\ }$$
(while keeping $r_+$, $\r_+$, $\xi_+$, $r_1$, $\r_1$, $\xi_1$ unvaried).
Indeed, since $1-e^{-x}\le x$ for all $x$,
$$\d_1=\d e^{-{\cal X}\|\frac{\o_\JJ}{\o_\rr}\|_{r,\r}} = \d-\d(1- e^{-{\cal X}\|\frac{\o_\JJ}{\o_\rr}\|_{r,\r}})\ge \d-{\cal X}\|\frac{\o_\JJ}{\o_\rr}\|_{r,\r}\ge \d-\d' =\d_{1\rm new}\ .$$
This also implies $\xi_+=\d_1-\d'\ge \d-2\d'=\xi_{+\rm new}$. That $s_+\ge s_{+\rm new}$, $s_1\ge s_{1\rm new}$ is even more immediate.\end{remark}

\nl
Now we can proceed with the

\paragraph{Proof of Lemma \ref{normal form lemma}}
Let $\widetilde{\tt c}_{n,m}$ be as in Lemma \ref{iterative lemma}. We shall choose $\ovl{\tt c}_{n,m}$ suitably large with respect to $\widetilde{\tt c}_{n,m}$.

\nl
We apply Lemma \ref{iterative lemma} with
$${2}r'=\frac{r}{3}\ ,\quad {2}\r'=\frac{\r}{3}\ ,\quad {2}\xi'=\frac{\xi}{3}\ ,\quad {3}s'=\frac{s}{3}\ ,\quad {3}\d'=\frac{\d}{3}\ ,\quad g\equiv0\ .$$
We make use of the stronger formulation described in Remark \ref{stronger iterative lemma}.  Conditions in \equ{ineq1} and the three former conditions in \equ{new cond} are trivially true. The two latter inequalities in 
\equ{new cond} reduce to 
$${\cal X}\|\frac{\o_\II}{\o_\rr}\|_{r,\r}<\frac{s}{9}\ ,\quad
{\cal X}\|\frac{\o_\JJ}{\o_\rr}\|_{r,\r}<
\frac{1}{9} $$
and they are certainly satisfied by assumption \equ{normal form assumptions}, for $N>1$. Since  
$$d=\min\{ \r' s', r'\xi', {\d'}^2 \}=\min\{ \r s/{36}, r\xi/{54}, {\d}^2/{81} \}\ge\frac{1}{{81}}
\min\{ \r s, r\xi, {\d}^2 \}=\frac{\eufm d}{{81}}
$$
we have that condition \equ{smallness} is certainly implied by the last inequality in \equ{normal form assumptions}, once one chooses ${\tt c}_{n,m}>{81} \widetilde{\tt c}_{n,m}$.
By Lemma \ref{iterative lemma}, it is then possible to conjugate $\HH$ to
$$\HH_1=\HH_0+\ovl f+f_1$$
with $f_1\in \cO_{r^\ppu,\r^\ppu,\xi^\ppu,s^\ppu,\d^\ppu}$, where $(r^\ppu,\r^\ppu,\xi^\ppu,s^\ppu,\d^\ppu):=2/3 (r, \r,\xi, s, \d)$ and
\beq{f1}\|f_1\|_{r^\ppu,\r^\ppu,\xi^\ppu,s^\ppu,\d^\ppu}\le {81}\widetilde{\tt c}_{n,m} \frac{{\cal X}}{\eufm d}   \|\frac{1}{\o_\rr}\|_{r,\r} \|\widetilde f\|_{r,\r,\xi,s,\d}\| f\|_{r,\r,\xi,s,\d}\le \frac{\| f\|_{r,\r,\xi,s,\d}}{2}\ .\eeq
since ${\tt c}_{n,m}\ge {162} \widetilde{\tt c}_{n,m}$ and $N\ge1$. Now we aim to apply Lemma \ref{iterative lemma} $N$ times, each time with parameters
$$r_j'=\frac{r}{6N}\ ,\quad \r_j'=\frac{\r}{6N}\ ,\quad \xi_j'=\frac{\xi}{6N}\ ,\quad s_j'=\frac{s}{{9}N}\ ,\quad \d_j'=\frac{\d}{{9}N}\ .$$
To this end, we let 
\beqano
&&r^{(j+1)}:=r^\ppu-j\frac{r}{3N}\ ,\quad \r^{(j+1)}:=\r^\ppu-j\frac{\r}{3N}\ ,\quad \xi^{(j+1)}:=\xi^\ppu-j\frac{\xi}{3N}\nonumber\\
&&s^{(j+1)}:=s^\ppu-j\frac{s}{3N}\ ,\quad \d^{(j+1)}:=\d^\ppu-j\frac{\d}{3N}\nonumber\\
&& r_1^\ppj:=r^\ppj\ ,\quad \r_1^\ppj:=\r^\ppj\ ,\quad \xi_1^\ppj:=\xi^\ppj\ ,\quad s_1^\ppj:=s^\ppj-\frac{s}{9N}\ ,\nonumber\\
&&\d_1^\ppj:=\d^\ppj-\frac{\d}{9N}\ ,\qquad {\cal X}_j:=\sup\{|x|:\ x\in \Xi_{\xi_j}\}
\eeqano
with $1\le j\le N$.

\nl
We assume that for a certain $1\le i\le N$ and all $1\le j\le i$, we have $\HH_j\in \cO_{r^\ppj, \r^\ppj, \xi^\ppj, s^\ppj, \d^\ppj}$ of the form
\beqa{Hi}
&&\HH_j=\HH_0+g_{j-1}+f_j\ , \quad g_{j-1}\in \cN_{r^\ppj, \r^\ppj, \xi^\ppj, s^\ppj, \d^\ppj}\ ,\quad g_{j-1}-g_{j-2}=\ovl f_{j-1}\\
&&\label{i+1 ineq} \|f_j\|_{r^\ppj, \r^\ppj, \xi^\ppj, s^\ppj, \d^\ppj}\le \frac{\|f_1\|_{r^\ppu, \r^\ppu, \xi^\ppu, s^\ppu, \d^\ppu}}{2^{j-1}}\eeqa
with  $g_{-1}\equiv0$,  $g_0=f_0=\ovl f$.
If $i=N$, we have nothing more to do. If $i<N$, we want to prove that
Lemma \ref{iterative lemma} can be applied so as to conjugate $\HH_i$ to a suitable $\HH_{i+1}$ such that \equ{Hi}--\equ{i+1 ineq}
are true with $j=i+1$.
To this end, we have to check
\beqa{smallness0}
&&{\cal X}_i\|\frac{\o_\II}{\o_\rr}\|_{r_i,\r_i}<s'_i
\ ,\quad
{\cal X}_i\|\frac{\o_\JJ}{\o_\rr}\|_{r_i,\r_i}<
\frac{\d'_i}{\d_i}\\
\label{smallness i}
&& \widetilde{\tt c}_{n,m}\frac{{\cal X}_i}{d_i}  \|\frac{1}{\o_\rr}\|_{r_i,\r_i}\|f_i\|_{r_i, \r_i, \xi_i, s_i, \d_i}<1\ .
\eeqa
where
$d_i:=\min\{ \r_i' s_i', r_i'\xi_i', {\d'}_i^2 \}$.
Conditions \equ{smallness0} are certainly  verified, since in fact they are implied by the definitions above 
(using also $\d_i\le \frac{2}{3}\d$, ${\cal X}_i\le {\cal X}$) and the two former inequalities in \equ{normal form assumptions}.
 To check the validity of \equ{smallness i}, we firstly observe that
$$d_i=\min\{r'_j\xi'_j,\ \r'_js'_j,\ (\d'_j)^2\}\ge \frac{\eufm d}{{81}N^2}\ .$$
Using then ${\tt c}_{n,m}>{162} \widetilde{\tt c}_{n,m}$,${\cal X}_i<{\cal X}$, Equation \equ{f1}, the inequality in \equ{i+1 ineq} with $j=i$ and the last inequality in \equ{normal form assumptions}, we easily conclude
\beqa{last but one}
&&\|f_i\|_{r_i, \r_i, \xi_i, s_i, \d_i}\le \|f_1\|_{r^\ppu, \r^\ppu, \xi^\ppu, s^\ppu, \d^\ppu}\le {81}\widetilde{\tt c}_{n,m} \frac{{\cal X}}{\eufm d}   \|\frac{1}{\o_\rr}\|_{r,\r} \| f\|^2_{r,\r,\xi,s,\d}
 \nonumber\\
 &&\le
 \frac{1}{\widetilde{\tt c}_{n,m}}\frac{\eufm d}{{81} N^2}\frac{1} {{\cal X}_i}  (\|\frac{1}{\o_\rr}\|_{r_i,\r_i})^{-1}\le  \frac{1}{\widetilde{\tt c}_{n,m}}\frac{d_i} {{\cal X}_i}  (\|\frac{1}{\o_\rr}\|_{r_i,\r_i})^{-1}\eeqa
 which is just \equ{smallness i}.

\nl
 Then the Iterative Lemma is applicable to $\HH_i$, and Equations \equ{Hi} with $j=i+1$ follow from it. The proof that also \equ{i+1 ineq} holds (for a possibly larger value of ${\tt c}_{n,m}$)  when $j=i+1$ proceeds along the same lines  as in \cite[proof of the Normal Form Lemma, p. 194--95]{poschel93} and therefore is omitted. The same for the proof of the first inequality in \equ{thesis}, for $g_N:=\HH_1$.
 \qed
 
%

\newpage

\def\cprime{$'$} \def\cprime{$'$} \def\cprime{$'$}

\end{document}